\documentclass{amsart}
\title{On Mackey decomposition for locally profinite groups}
\author{Yuki Yamamoto}
\usepackage{mathrsfs}
\usepackage{amssymb}
\usepackage[all]{xy}
\newtheorem{defn}{Definition}[section]
\newtheorem{rem}[defn]{Remark}
\newtheorem{thm}[defn]{Theorem}
\newtheorem{prop}[defn]{Proposition}
\newtheorem{lem}[defn]{Lemma}
\newtheorem{cor}[defn]{Corollary}
\newtheorem{eg}[defn]{Example}
\newtheorem{assump}[defn]{Assumption}

\newenvironment{prf}{Proof. }{\hfill$\Box$}

\DeclareMathOperator{\Hom}{Hom}
\DeclareMathOperator{\GL}{GL}

\DeclareMathOperator{\Aut}{Aut}

\DeclareMathOperator{\Res}{Res}
\DeclareMathOperator{\Ind}{Ind}
\DeclareMathOperator{\cInd}{c-Ind}

\DeclareMathOperator{\supp}{supp}

\newcommand{\ofra}{\mathfrak{o}}

\newcommand{\pfra}{\mathfrak{p}}

\newcommand{\N}{\mathbb{N}}
\newcommand{\Z}{\mathbb{Z}}

\newcommand{\C}{\mathbb{C}}

\begin{document}
\maketitle

\begin{abstract}
To study induced representation of some class of groups, Mackey's theory is very useful.  
In this paper, we consider some generalization of Mackey's theory for locally profinite groups.  
In particular, we give conditions on groups under which we have the Mackey decomposition and some examples such that we do not have the Mackey decomposition in some sense.  
\end{abstract}

\tableofcontents

\section{Introduction}

In representation theory for groups $G$, the induction $\Ind$ is fundamental as a method to obtain representations of $G$ from representations of its subgroups.  
To study induced representations, Mackey's theory is very useful.  
Mackey \cite{Ma} showed that for any finite group $G$ and its subgroups $H$ and $K$, we have the Mackey decomposition
\[
	\Res_{K}^{G} \Ind_{H}^{G} \rho \cong \bigoplus_{g \in K \backslash G / H} \Ind_{K \cap {}^{g}H}^{K} \Res_{K \cap {}^{g}H}^{{}^{g}H} {}^{g} \rho
\]
for any $H$-representation $\rho$ over some field, where for $g \in G$ we put ${}^{g}H:=gHg^{-1}$ and ${}^{g} \rho$ is a ${}^{g}H$-representation defined as ${}^{g} \rho (ghg^{-1})= \rho(h)$ for any $h \in H$.  
By Frobenius reciprocity, we can also examine intertwining operators between induced representations from the Mackey decomposition.  

There exist generalizations of Mackey's theory for some class of groups.  
The goal of this paper is to give a generalization of Mackey's theory for smooth representations of locally profinite groups, based on previous researches \cite{Ku}, \cite{Vi}.  

In smooth representation theory for locally profinite groups, there two functors treated as induction, ``induction'' $\Ind$ and ``compact induction'' $\cInd$.  
We give conditions on groups under which we have the Mackey decomposition for $\Ind$ and $\cInd$.  
The main theorem in this paper is the following.  

\begin{thm}[Theorem \ref{gen_of_MD}, Example \ref{ceg_for_gen_of_MD}]
\label{MD_for_intro}
Let $G$ be a locally profinite group.  
Let $H$ and $K$ be closed subgroups in $G$.  
Let $\rho$ be a smooth representation of $H$ over some commutative ring $R$ with unit.  
\begin{enumerate}
\item If either $H$ or $K$ is open in $G$, then we have the Mackey decomposition
\[
	\Res_{K}^{G} \cInd_{H}^{G} \rho \cong \bigoplus_{g \in K \backslash G / H} \cInd_{K \cap {}^{g}H}^{K} \Res_{K \cap {}^{g}H}^{{}^{g}H} {}^{g} \rho.  
\] \label{MD_for_cInd_for_intro}
\item If $K$ is open in $G$, then we have the Mackey decomposition 
\[
	\Res_{K}^{G} \Ind_{H}^{G} \rho \cong \left( \prod_{g \in K \backslash G / H} \Ind_{K \cap {}^{g}H}^{K} \Res_{K \cap {}^{g}H}^{{}^{g}H} {}^{g} \rho \right) ^{\infty}, 
\] \label{MD_for_Ind_for_intro}
where for a $K$-representation $\tau$, we denote by $\tau^{\infty}$ the $K$-smooth part of $\tau$.  
\item If we omit the assumption that $K$ is open, there exists an example such that $H$ is open in $G$ and the isomorphism in $(\ref{MD_for_Ind_for_intro})$ does not hold.  \label{ceg_for_MD_for_intro}
\end{enumerate}
\end{thm}

The example in (\ref{ceg_for_MD_for_intro}) shows that the condition that $HgK$ is open and closed in $G$ for any $g \in G$ is not enough to have the Mackey decomposition in general.  (cf. \cite[I.5.5]{Vi})

As a corollary of this theorem, we also consider intertwining operators between induced representations, which are studied by Kutzko \cite{Ku} for some special cases.  

\bigbreak

\noindent{\bfseries Acknowledgment}\quad
I am deeply grateful to my supervisor Naoki Imai for checking my draft carefully and giving helpful comments to improve the draft.  
I also thank Marie-France Vign\'eras and Masao Oi for many help, comments and suggestions.  
I am supported by the FMSP program at Graduate School of Mathematical Sciences, the University of Tokyo.  
I was also supported by JSPS KAKENHI Grant Number JP21J13751.  
\bigbreak

\section{Mackey decomposition}
\label{Mac_decomp}

In this section, we collect knowledge on the generalization of Mackey theory \cite{Ma} to locally profinite groups, based on \cite[I.5]{Vi}.  

\subsection{Induction and compact induction}
Let $G$ be a locally profinite group.  
Let $H$ be a closed subgroup in $G$, and let $\rho$ be a smooth $H$-representation over some (commutative) ring $R$ (with unit).  
Then we can consider the induction $\Ind_{H}^{G} \rho$ and the compact induction $\cInd_{H}^{G} \rho$ of $\rho$.  

We recall the definition of $\Ind_{H}^{G} \rho$ and $\cInd_{H}^{G} \rho$.  
Let $V$ be the representation space of $\rho$.  
We denote by $\mathrm{IND}_{H}^{G} (\rho)$ the set of functions $f \colon G \to V$ such that $f(hg)=\rho(h)f(g)$ for any $g \in G$ and $h \in H$.  
The set $\mathrm{IND}_{H}^{G} (\rho)$ can be equipped with a canonical $G$-representation structure by right translation.  
Let $\Ind_{H}^{G} \rho$ be the \textit{smooth} part of $\mathrm{IND}_{H}^{G} (\rho)$.  
Moreover, we also denote by $\cInd_{H}^{G} \rho$ the set of elements $f \in \Ind_{H}^{G} \rho$ such that $\supp f$ is compact modulo $H$.  
Then we can consider functors $\Ind_{H}^{G}$ and $\cInd_{H}^{G}$ from the category $\mathcal{M}_{R}(H)$ of smooth $H$-representations over $R$ to the category $\mathcal{M}_{R}(G)$ of smooth $G$-representations over $R$.  

Let $K$ be another closed subgroup in $G$.  
Then we can consider functors $\Res_{K}^{G} \Ind_{H}^{G}$ and $\Res_{K}^{G} \cInd_{H}^{G}$ from $\mathcal{M}_{R}(H)$ to $\mathcal{M}_{R}(K)$.  
We also can define other functors $\mathcal{M}_{R}(H) \to \mathcal{M}_{R}(K)$ as in \cite[I.5.4]{Vi} as we recall in the following.  
First, for $g \in G$ and an $H$-representation $(\rho,V)$, we can construct $K$-representations $\Ind_{K \cap {}^{g}H}^{K} \Res_{K \cap {}^{g}H}^{{}^{g}H} {}^{g}\rho$ and $\cInd_{K \cap {}^{g}H}^{K} \Res_{K \cap {}^{g}H}^{{}^{g}H} {}^{g}\rho$, where ${}^{g}H:=gHg^{-1}$ and ${}^{g} \rho$ is a ${}^{g}H$-representation defined as ${}^{g} \rho (h'):=\rho(g^{-1}h'g)$ for $h' \in {}^{g}H$.  
On the other hand, let $I(H,g,K; \rho)$ be the set of locally constant functions $f \colon HgK \to V$ such that $f(hx)=\rho(h)f(x)$ for $h \in H$ and $x \in HgK$.  
The set $I(H,g,K; \rho)$ can be equipped with a canonical $K$-representation structure by right translation.  
Let $I^{\infty}(H,g,K; \rho)$ be the smooth part of $I(H,g,K; \rho)$.  
Moreover, we also denote by $I_{c}^{\infty}(H,g,K; \rho)$ the set of elements $f \in I^{\infty}(H,g,K; \rho)$ such that $\supp f$ is compact modulo $H$.  
Then $\rho \mapsto I^{\infty}(H, g, K; \rho)$ and $\rho \mapsto I_{c}^{\infty}(H, g, K; \rho)$ are functors $\mathcal{M}_{R}(H) \to \mathcal{M}_{R}(K)$.  

The Mackey decomposition is a theorem which relates these functors when $H$ and $K$ are ``nice''.  
The following proposition is proved by arguments in \cite[I.5.4]{Vi}.  

\begin{prop}
\label{incl_of_MD}
Let $G$ be a locally profinite group, and let $H$ and $K$ be closed subgroups in $G$.  
Suppose $HgK$ is open in $G$ for any $g \in G$.  
Moreover, let $R$ be a commutative ring with unit, and let $\rho$ be an $H$-representation over $R$.  
\begin{enumerate}
\item Let $\phi_{H,g,K} \colon \mathrm{IND}_{H}^{G}(\rho) \to I(H,g,K; \rho)$ be a map defined by the restriction of functions $f \mapsto f |_{HgK}$.  
Then the map $\phi:=\prod_{g \in H \backslash G / K} \phi_{H,g,K} \colon \mathrm{IND}_{H}^{G}(\rho) \to \prod_{g \in H \backslash G / K} I(H,g,K; \rho)$ induces the following inclusions of $K$-representations:
\begin{enumerate}
\item $\Res_{K}^{G} \Ind_{H}^{G} \rho \hookrightarrow \prod_{g \in H \backslash G / K} I^{\infty}(H,g,K;\rho)$; \label{Ind_incl1}
\item $\Res_{K}^{G} \cInd_{H}^{G} \rho \hookrightarrow \bigoplus_{g \in H \backslash G / K} I_{c}^{\infty}(H,g,K; \rho)$.  \label{cInd_incl1}
\end{enumerate}
\item For $g \in G$, the map $I(H,g,K;\rho) \to \mathrm{IND}_{{}^{g^{-1}}H \cap K}^{K} \Res_{{}^{g^{-1}}H \cap K}^{{}^{g^{-1}}H} {}^{g^{-1}} \rho$ defined as $f \mapsto [ k \mapsto f(gk) ]$ induces an inclusion of $K$-representations 
\[
	I^{\infty}(H,g,K;\rho) \hookrightarrow \Ind_{K \cap {}^{g^{-1}}H}^{K} \Res_{K \cap {}^{g^{-1}}H} ^{{}^{g^{-1}}H} {}^{g^{-1}} \rho.  
\]
\label{Ind_isom0}
\item We have a canonical inclusion of $K$-representations 
\[
	\Res_{K}^{G} \Ind_{H}^{G} \rho \hookrightarrow \prod_{g \in K \backslash G / H} \Ind_{K \cap {}^{g}H}^{K} \Res_{K \cap {}^{g}H} ^{{}^{g}H} {}^{g} \rho.  
\]
\label{Ind_incl2}
\end{enumerate}
\end{prop}

Next, we consider these inclusions with some assumption on topological groups.  
We apply the assumption in the case as in Proposition $\ref{incl_of_MD}$.  
In particular, the set $HgK$ is open in $G$ for any $g \in G$.  

\begin{assump}
\label{assump_on_topgrp}
For any open subgroup $H'$ in $H$, any open subgroup $K'$ in $K$ and $g \in G$, the set $H'gK'$ is open in $G$.  
\end{assump}

\begin{prop}
\label{incl_of_MD2}
Let $G,H,K,R$ and $\rho$ be as in Proposition $\ref{incl_of_MD}$.  
Moreover, suppose that Assumption $\ref{assump_on_topgrp}$ holds.  
\begin{enumerate}
\item The inclusion in Proposition $\ref{incl_of_MD} (\ref{Ind_isom0})$ is isomorphic.  \label{Ind_isom1}
\item The isomorphism $(\ref{Ind_isom1})$ induces a $K$-isomorphism
\[
	I_{c}^{\infty}(H,g,K;\rho) \cong \cInd_{K \cap {}^{g^{-1}}H}^{K} \Res_{K \cap {}^{g^{-1}}H} ^{{}^{g^{-1}}H} {}^{g^{-1}} \rho.  
\]
\label{cInd_isom1}
\item The inclusion in Proposition $\ref{incl_of_MD} (\ref{Ind_incl2})$ induces a canonical inclusion of smooth $K$-representations
\[
\Res_{K}^{G} \cInd_{H}^{G} \rho \hookrightarrow \bigoplus_{g \in K \backslash G / H} \cInd_{K \cap {}^{g}H}^{K} \Res_{K \cap {}^{g}H} ^{{}^{g}H} {}^{g} \rho.  
\]
\label{cInd_incl2}
\end{enumerate}
\end{prop}

\begin{prf}
The claim (\ref{cInd_incl2}) follows from (\ref{cInd_isom1}) and Proposition \ref{incl_of_MD} (\ref{cInd_incl1}).  

We consider (\ref{Ind_isom1}).  
Let $V$ be the representation space for $\rho$.  
For a map $f \colon K \to V$ such that $f(hk)={}^{g^{-1}}\rho(h)f(k)$ for any $h \in g^{-1}Hg \cap K$ and $k \in K$, we can obtain a well-defined map $\hat{f} \colon HgK \to V$ as $\hat{f}(hgk)=\rho(h)f(k)$ for $h \in H$ and $k \in K$.  
The map $\hat{f}$ is the unique map $HgK \to V$ such that $\hat{f}(gk)=f(k)$ and $\hat{f}(hx)=\rho(h)\hat{f}(x)$ for any $h \in H, k \in K$ and $x \in HgK$.  
Moreover, if $f$ is $K$-smooth, then $\hat{f}$ is also $K$-smooth.  
Then, to show the inclusion in Proposition \ref{incl_of_MD} (\ref{Ind_isom0}) is surjective, it suffices to show that $\hat{f} \in I^{\infty}(H,g,K;\rho)$ for all $f \in \Ind_{K \cap {}^{g^{-1}}H}^{K} \Res_{K \cap {}^{g^{-1}}H} ^{{}^{g^{-1}}H} {}^{g^{-1}} \rho$, that is, the map $\hat{f}$ is locally constant.  

Let $x \in HgK$.  
Then we have $x=hgk$ for some $h \in H$ and $k \in K$, and so $\hat{f}(x)=\rho(h)f(k)$ by the definition of $\hat{f}$.  
Since $\rho$ is smooth, there exists an open subgroup $H'$ in $H$ that fixes $f(k)$.  
On the other hand, the map $f$ is invariant under the right translation by some open subgroup $K''$ in $K$.  
When we put $K'=kK''k^{-1}$, the subgroup $K'$ is also open in $K$.  
By Assumption \ref{assump_on_topgrp}, the set $H'gK'$ is open in $G$.  
Then $hH'gkK''=hH'gK'k$ is an open neighborhood of $hgk=x$ in $G$.  
Moreover, for $h' \in H'$ and $k'' \in K''$, we have $\hat{f}(hh'gkk'')=\rho(hh')f(kk'')=\rho(h)\rho(h')f(k)=\rho(h)f(k)=\hat{f}(x)$.  
Therefore, the map $\hat{f}$ is constant on $hH'gkK''$, and so $\hat{f}$ is locally constant.  

To show (\ref{cInd_isom1}), we consider the continuous bijection $\overline{\iota} \colon (g^{-1}Hg \cap K) \backslash K \to H \backslash HgK$ induced by $\iota \colon K \to HgK; \, k \mapsto gk$.  
By Assumption \ref{assump_on_topgrp}, for any open subgroup $K'$ in $K$, the set $HgK'$ is open in $G$.  
In particular, $H \backslash HgK'$ is open in $HgK$, which implies that $\overline{\iota}$ is an open map and so $\overline{\iota}$ is homeomorphic.  

Here, let $f \in \Ind_{K \cap {}^{g^{-1}}H}^{K} \Res_{K \cap {}^{g^{-1}}H} ^{{}^{g^{-1}}H} {}^{g^{-1}} \rho$.  
If we prove $Hg \supp f = \supp \hat{f}$ for any $f$, we have $\overline{\iota}\left( (g^{-1}Hg \cap K)\backslash \supp f \right) = H \backslash Hg\supp f = H \backslash \supp \hat{f}$.  
Since $\overline{\iota}$ is homeomorphic, the set $(g^{-1}Hg \cap K) \backslash K$ is compact if and only if $H \backslash \supp \hat{f}$ is compact, which implies ($\ref{cInd_isom1}$).  
Then, it suffices to show $Hg \supp f = \supp \hat{f}$.  

Let $x \in Hg \supp f$.  
Then we have $x=hgk$ for some $h \in H$ and $k \in \supp f$.  
Since $f$ is locally constant, we have $f(k) \neq 0$.  
Therefore $\hat{f}(x)=\rho(h)f(k) \neq 0$ since $\rho(h)$ is isomorphic, and so $x \in \supp \hat{f}$.  
On the other hand, let $x \in \supp \hat{f}$.  
Since $x \in \supp \hat{f} \subset HgK$, there exist $h \in H$ and $k \in K$ such that $x = hgk$.  
We have $\hat{f}(x) \neq 0$ since $\hat{f}$ is locally constant.  
Moreover, we also have $\hat{f}(x) = \hat{f}(hgk)=\rho(h)f(k)$ by definition of $\hat{f}$.  
Since $\rho(h)$ is isomorphic, we obtain $f(k) \neq 0$, and so $k \in \supp f$.  
Therefore we have $x \in Hg \supp f$, which completes the proof.  
\end{prf}

We consider sufficient conditions to satisfy the statement of Assumption \ref{assump_on_topgrp}.    

\begin{lem}
\label{lem_on_topgrp}
Let $G,H,K$ be as in Proposition $\ref{incl_of_MD}$.  
Moreover, suppose that $(G,H,K)$ satisfies either of the following:  
\begin{enumerate}
\item either $H$ or $K$ is open in $G$; 
\item the group $G$ is second-countable.  
\end{enumerate}
Then Assumption $\ref{assump_on_topgrp}$ holds.  
In particular, if $G$ is the group of $F$-rational points of some linear algebraic group over a non-archimedean local field $F$, then there exist the inclusions and isomorphisms as in Proposition $\ref{incl_of_MD2}$.  
\end{lem}

\begin{prf}
Suppose $H$ is open in $G$.  
Then any open subgroup $H'$ in $H$ is also open in $G$, and so $H'gK'$ is open in $G$ for any subset $K'$ in $K$ and $g \in G$.  
We can also prove the lemma for the case that $K$ is open in $G$ similarly.  

Suppose $G$ is second-countable.  
Then $H, K$ and $HgK$ are also second-countable.  
Let $H'$ (resp. $K'$) be an open subgroup in $H$ (resp. $K$).  
Since $H, K$ and $HgK$ are closed in $G$, the set $HgK$ is locally compact, Hausdorff and totally disconnected.  
Then there exists a compact open subgroup $H''$ (resp. $K''$) in $H'$ (resp. $K'$).  
If $H''gK''$ is open in $HgK$, then $H'gK'$, which is a union of open subsets in $G$ which are of the form $h'H''gK''k'$ for some $h' \in H'$ and $k' \in K'$, is also open in $HgK$.  
Since $HgK$ is open in $G$, the subset $H'gK'$ is also open in $G$.  
Therefore, we may assume $H'$ and $K'$ are compact.  

Since $H$ and $K$ are second-countable, we have $H=\bigcup_{i \in \N}h_{i}H'$ and $K=\bigcup_{j \in \N}K'k_{j}$ for some $h_{i} \in H$ and $k_{j} \in K$.  
Then $HgK = \bigcup_{i,j}h_{i}H'gK'k_{j}$.  
Since $H'gK'$ is compact, the subset $H'gK'$ is closed in $HgK$ as $HgK$ is Hausdorff.  
Therefore, by the Baire category theorem, there exist $(i,j) \in \N^{2}$ and a nonempty open subset $U$ in $HgK$ such that $U \subset h_{i}H'gK'k_{j}$.  
Then $U':=h_{i}^{-1}Uk_{j}^{-1}$ is a nonempty subset in $H'gK'$ which is open in $G$.  
Since $H'gK'=H'U'K'$, the set $H'gK'$ is open in $HgK$ and $G$.  
\end{prf}

\subsection{Mackey decomposition and some counterexamples}
In general, the inclusion $(\ref{Ind_incl2})$ in Proposition $\ref{incl_of_MD}$ is not necessarily isomorphic, as shown later.  
Our main proposition gives us some sufficient conditions to make the inclusions in Proposition $\ref{incl_of_MD}$ and Proposition $\ref{incl_of_MD2}$ isomorphic.  

\begin{prop}[Mackey decomposition]
\label{gen_of_MD}
Let $G, H, K, R$ and $\rho$ be as in Proposition $\ref{incl_of_MD}$.  
\begin{enumerate}
\item Suppose $K$ is open.  
Then we have 
\[
\Res_{K}^{G} \Ind_{H}^{G} \rho \cong \left( \prod_{g \in K \backslash G / H} \Ind_{K \cap {}^{g}H}^{K} \Res_{K \cap {}^{g}H} ^{{}^{g}H} {}^{g} \rho \right)^{\infty}, 
\]
where for a $K$-representation $\tau$, the representation $\tau^{\infty}$ is the $K$-smooth part of $\tau$.  \label{gen_of_MD_for_Ind}
\item Suppose either $H$ or $K$ is open in $G$.  
Then we have the Mackey decomposition 
\[
\Res_{K}^{G} \cInd_{H}^{G} \rho \cong \bigoplus_{g \in K \backslash G / H} \cInd_{K \cap {}^{g}H}^{K} \Res_{K \cap {}^{g}H} ^{{}^{g}H} {}^{g} \rho.  
\] \label{gen_of_MD_for_cInd}
\end{enumerate}
\end{prop}

\begin{prf}
We show (\ref{gen_of_MD_for_Ind}).  
Since $G$-smooth representations are also $K$-smooth, the inclusion (\ref{Ind_incl1}) in Proposition \ref{incl_of_MD} induces an inclusion
\[
	\Res_{K}^{G} \Ind_{H}^{G} \rho \hookrightarrow \left( \prod_{g \in H \backslash G / K}  I^{\infty}(H,g,K;\rho) \right)^{\infty}.  
\]
By (\ref{Ind_isom1}) in Proposition \ref{incl_of_MD2} and Lemma \ref{lem_on_topgrp}, it is enough to show that the above inclusion is surjective.  
Let $(f_{g})_{g \in H \backslash G / K} \in  \left( \prod_{g \in H \backslash G / K}  I^{\infty}(H,g,K;\rho) \right)^{\infty}$.  
Then there exists an open subgroup $K'$ in $K$ such that $K'$ fixes $f_{g}$ for all $g \in H \backslash G / K$.  
We define a map $f \colon G \to V$ as $f(x)=f_{g}(x)$ for $x \in HgK$.  
Then $f \in \mathrm{IND}_{H}^{G}\rho$ and $\phi(f) = (f_{g})_{g \in H \backslash G / K}$.  
Therefore, it is enough to show that $f \in \Ind_{H}^{G} \rho$, that is, the element $f$ is $G$-smooth.  

Since $K$ is open in $G$, the open subgroup $K'$ in $K$ is also open in $G$.  
We show that $K'$ fixes $f$, which implies $f$ is $G$-smooth.  
Let $x \in G$.  
Then $x \in HgK$ for some double coset $HgK$.  
For $k' \in K'$, we have $f(xk')=f_{g}(xk')$ as $xk' \in HgKk'=HgK$.  
Since $f_{g}$ is fixed by $K'$, we also have $f_{g}(xk')=f_{g}(x)=f(x)$.  

We show (\ref{gen_of_MD_for_cInd}).  
It is enough to show that the inclusion (\ref{cInd_incl1}) in Proposition \ref{incl_of_MD} is surjective.  
Let $g \in G$, and let $f \in I_{c}^{\infty}(H,g,K; \rho)$.  
Let $V$ be the representation space of $\rho$.  
We define a map $\tilde{f} \colon G \to V$ as
\[
	\tilde{f}(x)= \left\{ 
	\begin{array}{cl}
	f(x) & (x \in HgK) \\
	0 & (x \notin HgK).  
	\end{array}
	\right.
\]
Then $\tilde{f} \in \mathrm{IND}_{H}^{G} \rho$ and $\phi(\tilde{f})=f$.  
Moreover, we have $\supp \tilde{f} = \supp f$, where the latter is compact modulo $H$.  
Therefore, it suffices to show that $\tilde{f}$ is a smooth element in $\mathrm{IND}_{H}^{G} \rho$, that is, there exists an open subgroup $J$ in $G$ such that $\tilde{f}$ is $J$-invariant.  

Suppose $K$ is open.  
Then any open subgroup $K'$ in $K$ which fixes $f$ is also open in $G$.  
By the assumption of $K'$, we have $\tilde{f}(xk')=f(xk')=f(x)=\tilde{f}(x)$ for any $x \in HgK$ and $k' \in K'$.  
Moreover, for this $K'$, we have $(G \setminus HgK) \cdot K' = G \setminus HgK$, which implies $\tilde{f}(xk')=0=\tilde{f}(x)$ for any $x \in G \setminus HgK$ and $k' \in K'$.  
Therefore it is enough to put $J=K'$.  

Suppose $H$ is open.  
Then the set $H \backslash \supp \tilde{f}=H \backslash \supp f$ is discrete.  
Since $H \backslash \supp f$ is also compact, the set $H \backslash \supp \tilde{f}$ is finite.  
Therefore, there exist finitely many elements $x_{1}, \ldots, x_{n} \in G$ such that $\supp \tilde{f}$ is the disjoint union $\bigcup_{i=1}^{n} Hx_{i}$.  
For any $i$, there exists a compact open subgroup $H_{i}$ in $H$ which fixes $f(x_{i})$.  
We put $J:= \bigcap_{i=1}^{n}x_{i}^{-1}H_{i}x_{i}$, which is an open subgroup in $G$.  

We will show $\tilde{f}$ is $J$-invariant.  
First, we will show $\supp\tilde{f} \cdot j = \supp\tilde{f}$ for any $j \in J$.  
We have $\supp\tilde{f} \cdot j = \bigcup_{i=1}^{n}Hx_{i}j$.  
For any $i$, there exists $h'_{i} \in H_{i}$ such that $j = x_{i}^{-1}h_{i}x_{i}$ since $J = \bigcap_{i=1}^{n}x_{i}^{-1}H_{i}x_{i}$.  
Then we obtain $\bigcup_{i=1}^{n}Hx_{i}j = \bigcup_{i}Hx_{i} \cdot x_{i}^{-1}h_{i}x_{i} = \bigcup_{i}Hx_{i} = \supp\tilde{f}$.  
In particular, we also have $(G \setminus \supp \tilde{f}) \cdot J = G \setminus \supp \tilde{f}$, which implies that $\tilde{f}(xj)=0=\tilde{f}(x)$ for any $x \in G \setminus \supp \tilde{f}$ and $j \in J$.  
On the other hand, let $x \in \supp \tilde{f}$ and $j \in J$.  
Then there exists $i$ such that $x=hx_{i}$ for some $h \in H$ and $j=x_{i}^{-1}h_{i}x_{i}$ for some $h_{i} \in H_{i}$.  
Therefore, we have $\tilde{f}(xj) = \tilde{f}(hh_{i}x_{i}) = \rho(hh_{i})\tilde{f}(x_{i})$.  
Here, we have $\rho(hh_{i})\tilde{f}(x_{i})=\rho(h)\left( \rho(h_{i})f(x_{i}) \right)=\rho(h)f(x_{i})$ since $h_{i} \in H_{i}$ fixes $f(x_{i})$.  
Then we obtain $\tilde{f}(xj) = \rho(h)f(x_{i})=f(hx_{i})=f(x)=\tilde{f}(x)$, which completes the proof.  
\end{prf}

In the following, we construct some example such that the inclusion (\ref{Ind_incl1}) in Proposition \ref{incl_of_MD} is not surjective.  

\begin{eg}
\label{ceg_for_gen_of_MD}
Let $F$ be a non-archimedean local field, and let $E/F$ be a ramified quadratic extension of fields.  
We denote by $\ofra_{F}$ (resp. $\ofra_{E}$) the ring of integers in $F$ (resp. $E$.)
We also denote by $\pfra_{F}$ (resp. $\pfra_{E}$) the maximal ideal in $\ofra_{F}$ (resp. $\ofra_{E}$.)
We put $G=\GL_{2}(E)$ and $K=\GL_{2}(F)$.  
We also put $H= \left( \begin{array}{cc} 1+\pfra_{E} & \pfra_{E} \\ \pfra_{E} & 1+\pfra_{E} \end{array} \right)$, which is a compact open subgroup in $G$.  
Let $(\rho, R)$ be the trivial representation of $H$ for some commutative ring $R$ with unit.  
Then $G, H, K, R$ and $\rho$ satisfies the condition in Proposition $\ref{incl_of_MD}$.  

We define $f \colon HK \to R$ as $f(x)=1$ for $x \in HK$.  
Then $f \in I^{\infty}(H, 1_{G}, K; \rho) \subset \left( \prod_{g \in H \backslash G/K} I^{\infty}(H, g, K; \rho) \right)^{\infty}$, and we can define $\tilde{f} \in \mathrm{IND}_{H}^{G} \rho$ such that $\phi(\tilde{f})=f$ as in the proof of Proposition $\ref{gen_of_MD}$.  
Here, by the definition of $\tilde{f}$, the function $\tilde{f}$ is the characteristic function of $HK$ in $G$.  

In the following, we will show that $\tilde{f} \notin \Ind_{H}^{G} \rho$, that is, the function $\tilde{f}$ is not fixed by any open subgroup $U $ in $G$.  
Since $\tilde{f}$ is the characteristic function of $HK$ in $G$, it is enough to show that $HKU \neq HK$ for any open subgroup $U$ in $G$.  
Let $U$ be an open subgroup in $G$.  
We put $u_{n}= \left( \begin{array}{cc} 1 & \varpi_{E}^{2n+1} \\ 0 & 1 \end{array} \right)$ for $n \in \N$, where we fix a uniformizer $\varpi_{E}$ of $E$.  
Since $\lim_{n \to \infty} u_{n} = 1_{G}$, we have $u_{n} \in U$ for some $n \in \N$.  
Then $\left( \begin{array}{cc} 1 & 0 \\ 0 & \varpi_{F}^{n+1} \end{array} \right) u_{n} \in KU \subset HKU$, where we also fix a uniformizer $\varpi_{F}$ of $F$.  
Here, suppose that $HKU=HK$.  
Then we have $\left( \begin{array}{cc} 1 & \varpi_{F}^{-(n+1)} \varpi_{E}^{2n+1} \\ 0 & 1 \end{array} \right) = \left( \begin{array}{cc} 1 & 0 \\ 0 & \varpi_{F}^{n+1} \end{array} \right) u_{n} \cdot \left( \begin{array}{cc} 1 & 0 \\ 0 & \varpi_{F}^{-(n+1)} \end{array} \right) \in HKU \cdot K = HK$.  
Therefore, there exists $h=\left( \begin{array}{cc} 1+a & b \\ c & 1+d \end{array} \right) \in H$ such that $h \left( \begin{array}{cc} 1 & \varpi_{F}^{-(n+1)} \varpi_{E}^{2n+1} \\ 0 & 1 \end{array} \right) \in K$.  
However, we have $h \left( \begin{array}{cc} 1 & \varpi_{F}^{-(n+1)} \varpi_{E}^{2n+1} \\ 0 & 1 \end{array} \right) = \left( \begin{array}{cc} x' & (1+a)\varpi_{F}^{-(n+1)} \varpi_{E}^{2n+1} +b \\ z' & w' \end{array} \right)$ for some $x', z', w' \in E$.  
When we define $v_{E}$ as the normalized valuation $E \to \Z \cup \{ \infty \}$, we have $v_{E}\left( (1+a)\varpi_{F}^{-(n+1)} \varpi_{E}^{2n+1} +b \right) =-1$ since $E/F$ is ramified quadratic and $a, b \in \pfra_{E}$.  
Since $-1 \notin 2 \Z \cup \{ \infty \} = v_{E}(F)$, we have $(1+a)\varpi_{F}^{-(n+1)} \varpi_{E}^{2n+1} +b \notin F$, which is a contradiction.  
\end{eg}

The Mackey decomposition leads to the following corollary on intertwining operators.  

\begin{cor}
\label{cor_of_gen_of_MD}
Let $G$ be a locally profinite group, $H$ be a closed subgroup in $G$, and let $K$ be an open subgroup in $G$.  
Let $\sigma$ (resp. $\tau$) be a smooth representation of $K$ (resp. $H$) over some commutative ring $R$ with unit.  
\begin{enumerate}
\item We have $\Hom_{G}(\cInd_{K}^{G} \sigma, \Ind_{H}^{G} \tau) \cong \prod_{g \in K \backslash G / H} \Hom_{K \cap {}^{g}H}(\sigma, {}^{g}\tau)$.  \label{cor_of_gen_of_MD_for_Ind}
\item Moreover, we further assume that $H \backslash HgK$ is compact for any $g \in G$ and $\sigma$ is finitely generated as a $K$-representation.  
Then we also have $\Hom_{G}(\cInd_{K}^{G} \sigma, \cInd_{H}^{G} \tau) \cong \bigoplus_{g \in K \backslash G / H} \Hom_{K \cap {}^{g}H} (\sigma, {}^{g} \tau)$.  \label{cor_of_gen_of_MD_for_cInd}
\end{enumerate}
\end{cor}

\begin{prf}
We show (\ref{cor_of_gen_of_MD_for_Ind}).  
Since $K$ is open in $G$, we have $\Hom_{G}(\cInd_{K}^{G} \sigma, \Ind_{H}^{G} \tau) \cong \Hom_{K}( \sigma, \Res_{K}^{G} \Ind_{H}^{G} \tau)$ by Frobenius reciprocity for compact induction.  
By Proposition \ref{gen_of_MD} (\ref{gen_of_MD_for_Ind}), we also have $\Res_{K}^{G} \Ind_{H}^{G} \tau \cong \left( \prod_{g \in K \backslash G / H} \Ind_{K \cap {}^{g}H}^{K} \Res_{K \cap {}^{g}H}^{{}^{g}H} {}^{g} \tau \right)^{\infty}$.  
Since $\sigma$ is smooth, we obtain 
\begin{eqnarray*}
\Hom_{K}(\sigma, \Res_{K}^{G} \Ind_{H}^{G} \tau) & \cong & \Hom_{K} \left( \sigma,  \left( \prod_{g \in K \backslash G / H} \Ind_{K \cap {}^{g}H}^{K} \Res_{K \cap {}^{g}H}^{{}^{g}H} {}^{g} \tau \right)^{\infty} \right) \\
 & \cong & \Hom_{K} \left( \sigma,  \prod_{g \in K \backslash G / H} \Ind_{K \cap {}^{g}H}^{K} \Res_{K \cap {}^{g}H}^{{}^{g}H} {}^{g} \tau \right) \\ 
 & \cong & \prod_{g \in K \backslash G / H} \Hom_{K} (\sigma, \Ind_{K \cap {}^{g}H}^{K} \Res_{K \cap {}^{g}H}^{{}^{g}H} {}^{g} \tau).  
\end{eqnarray*}
By Frobenius reciprocity for induction, we obtain 
\[
\Hom_{K}(\sigma, \Ind_{K \cap {}^{g}H}^{K} \Res_{K \cap {}^{g}H}^{{}^{g}H} {}^{g} \tau) \cong \Hom_{K \cap {}^{g}H} (\sigma, {}^{g} \tau)
\]
and complete the proof.    

We show (\ref{cor_of_gen_of_MD_for_cInd}).  
We have $\Hom_{G}(\cInd_{K}^{G} \sigma, \cInd_{H}^{G} \tau) \cong \Hom_{K}( \sigma, \Res_{K}^{G} \cInd_{H}^{G} \tau)$ as (\ref{cor_of_gen_of_MD_for_Ind}).  
Here, we also have $\Res_{K}^{G} \cInd_{H}^{G} \tau \cong \bigoplus_{g \in K \backslash G / H} \cInd_{K \cap {}^{g}H}^{K} \Res_{K \cap {}^{g}H}^{{}^{g}H} {}^{g} \tau$ by Proposition \ref{gen_of_MD} (\ref{gen_of_MD_for_cInd}).  
Since $\sigma$ is finitely generated, we obtain 
\begin{eqnarray*}
\Hom_{G}(\cInd_{K}^{G} \sigma, \cInd_{H}^{G} \tau) & \cong & \Hom_{K} \left( \sigma, \bigoplus_{g \in K \backslash G / H} \cInd_{K \cap {}^{g}H}^{K} \Res_{K \cap {}^{g}H}^{{}^{g}H} {}^{g} \tau \right) \\
 & \cong & \bigoplus_{g \in K \backslash G / H} \Hom_{K} \left( \sigma, \cInd_{K \cap {}^{g}H}^{K} \Res_{K \cap {}^{g}H}^{{}^{g}H} {}^{g} \tau) \right).  
\end{eqnarray*}
Here, since $K$ is open, the continuous bijection $(K \cap {}^{g}H) \backslash K \to gHg^{-1} \backslash gHg^{-1}K \cong H \backslash Hg^{-1}K$ is homeomorphic.  
Then, by the assumption that $H \backslash Hg^{-1}K$ is compact, the set $(K \cap {}^{g}H) \backslash K$ is also compact, and so the functor $\cInd_{K \cap {}^{g}H}^{K}$ is equal to $\Ind_{K \cap {}^{g}H}^{K}$.  
Therefore we obtain
\begin{eqnarray*}
\Hom_{K} \left( \sigma, \cInd_{K \cap {}^{g}H}^{K} \Res_{K \cap {}^{g}H}^{{}^{g}H} {}^{g} \tau \right)  & = & \Hom_{K} \left( \sigma, \Ind_{K \cap {}^{g}H}^{K} \Res_{K \cap {}^{g}H}^{{}^{g}H} {}^{g} \tau \right) \\
 & \cong & \Hom_{K \cap {}^{g}H} (\sigma, {}^{g} \tau), 
\end{eqnarray*}
where the last equation follows from Frobenius reciprocity for induction.  
\end{prf}

\begin{rem}
\begin{enumerate}
\item When $R=\C$ and $\sigma$ is a finite-dimensional representation of $K$, this corollary is proved by Kutzko \cite[Theorem 1.1]{Ku}.  
\item In the statement of \cite[Theorem 1.1]{Ku}, line 8, the word ``$I(\sigma_{c}^{G},\tau^{G})$'' seems to be a typographical error and seems to need to be replaced with ``$I(\sigma_{c}^{G},\tau_{c}^{G})$'' based on the proof.  
Even though this typo is corrected, the assumption that the image of $K$ in $H \backslash G$ is compact is not sufficient to make the theorem hold.  
The following example is a counterexample for the theorem.  
\end{enumerate}
\end{rem}

\begin{eg}
Let $\varphi \colon \Z/2\Z \to \Aut(\Z^{2})$ be the group homomorphism defined as $\varphi(1+2\Z)(n_{1}, n_{2})=(n_{2}, n_{1})$ for $(n_{1}, n_{2}) \in \Z^{2}$.  
We put $G=\Z^{2} \rtimes_{\varphi} (\Z/2\Z)$ with the discrete topology.  
In particular, $G$ is a locally profinite group.  
We put $H=K=(\Z \times \{ 0 \}) \rtimes \{ 0 \}$.  
We also put $g_{n}:=\left( (0,n),0 \right) \in G$ for $n \in \Z$ and $g':= \left( (0,0), 1+2\Z \right) \in G$.  
Moreover, let $\chi_{1}, \, \chi_{2}$ be distinct characters $H \to \C^{\times}$.  

We check that the tuple $(G,H,K,R=\C,\sigma=\chi_{1},\tau=\chi_{2})$ satisfies the condition which is assumed in Kutzko's statement.  
First, $H$ is a closed subgroup in $G$ and $K=H$ contains a compact open subgroup $\{ 1 \}$ in $G$ since $G$ is discrete.  
Next, $\sigma=\chi_{1}$ and $\tau=\chi_{2}$ are smooth, finite-dimensional representations of $K=H$.  
Moreover, the image of $K$ to $H \backslash G$ is $H \backslash HK=H \backslash H$, which is compact.  

Here, we consider the set of double cosets $K \backslash G / H$.  
For $n \in \Z$, we have $Kg_{n}H=( \Z \times \{ n \} ) \rtimes \{ 0 \}$.  
On the other hand, we also have $Kg'H=\Z^{2} \rtimes \{ 1+2\Z \}$.  
Therefore, the set $G$ is the disjoint union of $Kg_{n}H$ for $n \in \Z$ and $Kg'H$, that is, the set $\{ g_{n} \mid n \in \Z \} \cup \{ g' \}$ is a representative of $K \backslash G / H$.  

We show $\Hom_{G}(\cInd_{K}^{G} \sigma, \cInd_{H}^{G} \tau) \not\cong \bigoplus_{g \in K \backslash G / H} \Hom_{K \cap {}^{g}H} (\sigma, {}^{g} \tau)$.  
On the left-hand side, we have 
\[
\Hom_{G}(\cInd_{K}^{G} \sigma, \cInd_{H}^{G} \tau) \cong \bigoplus_{g \in K \backslash G / H} \Hom_{K}(\sigma, \cInd_{K \cap {}^{g}H}^{K} \Res_{K \cap {}^{g}H}^{{}^{g}H} {}^{g} \tau)
\]
as in the proof of Corollary $\ref{cor_of_gen_of_MD} (\ref{cor_of_gen_of_MD_for_cInd})$.  
If $g \in Kg_{n}H$, then $g$ centralizes $H$, and so ${}^{g}H=H$ and ${}^{g} \tau = \tau$.  
Then we have $\Hom_{K}(\sigma, \cInd_{K \cap {}^{g}H}^{K} \Res_{K \cap {}^{g}H}^{{}^{g}H} {}^{g} \tau)=\Hom_{K}(\sigma, \tau)=0$.  
If $g \in Kg'H$, then ${}^{g}H=( \{ 0 \} \times \Z) \rtimes \{ 0 \}$, and so $K \cap {}^{g}H= \{ 1_{G} \}$.  
Therefore we have $\cInd_{K \cap {}^{g}H}^{K} \Res_{K \cap {}^{g}H}^{{}^{g}H} {}^{g} \tau = \cInd_{ \{1_{G} \} }^{K} \mathbf{1} \cong \cInd_{ \{ * \} }^{\Z} \mathbf{1} = \{ f \colon \Z \to \C \mid | \supp f | < \infty \}$.  
Then to show $\Hom_{K} (\sigma, \cInd_{K \cap {}^{g}H}^{K} \Res_{K \cap {}^{g}H}^{{}^{g}H} {}^{g} \tau) = 0$, it suffices to show that $\cInd_{ \{ * \} }^{\Z} \mathbf{1}$ does not have any nontrivial, finite-dimensional subrepresentation of $\Z$.  
Suppose $W$ is a nonzero, finite-dimensional subspace in $\cInd_{ \{ * \} }^{\Z} \mathbf{1}$.  
Then the minimum $n_{0}:= \min \{ n \in \Z \mid f(n) \neq 0 \text{ for some } f \in W \}$ exists.  
Let $f \in W$ such that $f(n_{0}) \neq 0$.  
Then, we have $\left( \left( \cInd_{ \{ * \} }^{\Z} \mathbf{1} (1) \right) f \right) (n_{0}-1) = f \left( (n_{0}-1)+1 \right)=f(n_{0}) \neq 0$.  
By the definition of $n_{0}$, we obtain $\left( \cInd_{ \{ * \} }^{\Z} \mathbf{1} (1) \right) f \notin W$.  
In particular, $W$ is not invariant by the action of $\Z$, that is, $\cInd_{ \{ * \} }^{\Z} \mathbf{1}$ does not have any nontrivial, finite-dimensional subrepresentations of $\Z$.  
Therefore we obtain $\Hom_{G}(\cInd_{K}^{G} \sigma, \cInd_{H}^{G} \tau) = 0$.  

On the right-hand side, let $g \in Kg_{n}H$ for some $n \in \Z$.  
Then, as the left-hand side, we have $K \cap {}^{g}H = H, \, {}^{g} \tau = \tau$ and so $\Hom_{K \cap {}^{g}H}(\sigma, {}^{g} \tau) = \Hom_{H}(\sigma, \tau) = 0$.  
On the other hand, suppose that $g \in Kg'H$.  
Then $K \cap {}^{g}H = \{ 1_{G} \}$, and so $\Hom_{K \cap {}^{g}H} (\sigma, {}^{g} \tau) = \Hom_{\C}(\C, \C) \cong \C$.  
Therefore, $\bigoplus_{g \in K \backslash G / H} \Hom_{K \cap {}^{g}H} (\sigma, {}^{g} \tau) \cong \C \neq 0$.  
\end{eg}

\bigbreak\bigbreak
\noindent Yuki Yamamoto\par
\noindent Graduate School of Mathematical Sciences, 
The University of Tokyo, 3--8--1 Komaba, Meguro-ku, 
Tokyo, 153--8914, Japan\par
\noindent E-mail address: \texttt{yukiymmt@ms.u-tokyo.ac.jp}

\end{document}